\documentclass[a4paper,10,5pt,final,leqno,notitlepage]{article}

\usepackage{amsmath}
\usepackage{amsfonts}
\usepackage{amssymb}
\usepackage{amsthm} 

\usepackage[blocks]{authblk}

\usepackage[T1]{fontenc}

\usepackage{fullpage}

\usepackage{enumerate}
\usepackage{enumitem}

\begin{document}
\renewcommand{\thepage}{\small\arabic{page}}
\renewcommand{\thefootnote}{(\arabic{footnote})}
\renewcommand{\thesection}{\arabic{chapter}.\arabic{section}}
\renewcommand{\thesubsection}{\arabic{subsection}}

\renewcommand\Affilfont{\small}

\author{Maria Trybu\l{}a}
\affil{Institute of Mathematics, Faculty of Mathematics and Computer Science,\\ Jagiellonian University, \L ojasiewicza 6, 30-348 Krak\'ow, Poland\\ maria.trybula@im.uj.edu.pl}
\title{Proper holomorphic mappings, Bell's formula and the Lu Qi-Keng
problem on tetrablock}
\maketitle

\begin{abstract}
We consider a proper holomorphic map $\pi:D\rightarrow G$ between domains in $\mathbb{C}^{n}$ and show that it induces a unitary isomorphism between the Bergman space $\mathbb{A}^{2}(G)$ and some subspace of $\mathbb{A}^{2}(D)$. Using this isomorphism we construct orthogonal projection onto that subspace and we derive Bell's transformation formula for the Bergman kernel under proper holomorphic mappings. As a consequence of the formula we get that the tetrablock is not a Lu Qi-Keng domain.

\end{abstract}

\newtheorem{lemma}{Lemma}
\newtheorem{propozycja}{Proposition}
\newtheorem{corollary}{Corollary}
\newtheorem{Theorem}{Theorem}
\newtheorem{uwaga}{Remark}
\subsection{Introduction}

In this paper we consider Bergman Space $\mathbb{A}^{2}_{\phi}(G)$ consisting of  square integrable holomorphic functions on domain $G$ in $\mathbb{C}^{n}$ with weight $\phi$, that is $\int_{G}|f|^{2}\phi\, dV<\infty$,  where $\phi\in {L}^{\infty}_{\textup{loc}}(G)$ is some positive function. Let $D$ be another domain in $\mathbb{C}^{n}$ and suppose there is some   
proper holomorphic map $\pi:D\rightarrow G$. Let $J\pi$ denote the complex Jacobian of $\pi$. We show that $\mathbb{A}_{\phi\circ \pi}^{2}(D)$ has some closed subspace $H$ induced by $\pi$ which is unitary isomorphic to $\mathbb{A}_{\phi}^{2}(G)$. We also derive the formula for orthogonal projection onto $H$. Briefly, we construct a unitary operator $\Gamma:\mathbb{A}_{\phi}^{2}(G)\rightarrow\mathbb{A}_{\phi\circ \pi}^{2}(D)$ which is defined as follows $\Gamma f=\frac{1}{\sqrt{m}}(f\circ\pi)J\pi$, where $m$ stands for the multiplicity of $\pi$. In the sequel we use adjoint operator to $\Gamma$, which is in fact equal to $\Gamma^{-1}$, if $\Gamma$ is understood as an operator from $\mathbb{A}^{2}_{\phi}(G)$ to $\Gamma\mathbb{A}^{2}_{\phi}(G)$, therefore we show how $\Gamma^{*}$ acts on $\Gamma\mathbb{A}^{2}_{\phi}(G)$. More precisely If we take any $g\in\Gamma \mathbb{A}_{\phi}^{2}(G)$ then $\frac{g}{J\pi}$ is a well-defined function on a dense open subset of $D$ (the set of regular values of $\pi$), which is additionally invariant under $\pi$, that is if $z,\,w\in D$ are such that $\pi(z)=\pi(w),\,J\pi(w),J\pi(z)\not= 0$, then $\frac{g}{J\pi}(\pi(z))=\frac{g}{J\pi}(\pi(w))$. Therefore, equality $\widetilde{(\frac{g}{J\pi})}(z)=\frac{g}{J\pi}(\pi(z))$, where $z\in D$, defines well a holomorphic function on $D$ except for the (analytic) set of critical points of $\pi$. But the Riemann Removable Singularity Theorem for square integrable holomorphic functions ensures that $\widetilde{(\frac{g}{J\pi})}$ has a holomorphic extension on $D$. Having this, adjoint operator to $\Gamma$ might be described by equation $\Gamma^{*}g=\sqrt{m}\widetilde{(\frac{g}{J\pi})}$.  And now more precisely:
\begin{Theorem}
Let $\pi:D\rightarrow G$ be a proper holomorphic map between domains  in $\mathbb{C}^{n}$ with multiplicity m, $\phi\in L^{\infty}_{loc}(G,\mathbb{R}_{>0})$. Then the Hilbert subspace $\Gamma\mathbb{A}_{\phi}^{2}(G)$ of $\mathbb{A}_{\phi\circ \pi}^{2}(D)$ is isometrically isomorphic with $\mathbb{A}_{\phi}^{2}(G)$. The orthogonal projection P onto $\Gamma\mathbb{A}_{\phi}^{2}(G)$ is given by the formula
$$Pg=\frac{1}{m}\sum_{k=1}^{m}(g\circ \pi^{k}\circ\pi)J(\pi^{k}\circ \pi)$$
for $g\in\mathbb{A}^{2}(D)$, where $\{\pi^{j}\}_{j=1}^{m}$ are the local inverses to $\pi$.
\end{Theorem}
Note that it will follow from the proof that the formula on the right side actually defines a function from $\Gamma\mathbb{A}_{\phi}^2(G)\subset \mathbb{A}_{\pi\circ\phi}^2(D)$.

Recall that the Bergman Kerenel $K_{\phi}$ with weight $\phi$  is the reproducing kernel for the space $\mathbb{A}^{2}_{\phi}(G)$ that is all functionals of evaluation $ev_{z}:f\rightarrow f(z)$ are bounded for every $z\in G$ so from Riesz Representation Theorem there is unique $K_{z}\in\mathbb{A}^{2}_{\phi}(G)$ s.t $ev_{z}(f)=\langle f,K_{z}\rangle_{\mathbb{A}^{2}_{\phi}(G)}$, and $K_{\phi}(z,w)=\langle K_{w},K_{z}\rangle_{\mathbb{A}^{2}_{\phi}(G)}$. For $\phi=1$ we simply write $K_{\phi}=K$.
Definition and basic properties of the Bergman kernel might be found for instance in \cite{Jarnicki}. As a corollary of Theorem 1 we get Bell's Theorem. 
\begin{corollary}[see \cite{Bell}]
Let D, G be domains in $\mathbb{C}^{n}$ and let $\pi :D\rightarrow G$ be a proper holomorphic map with multiplicity $m$. Denote by $\pi^{1},\ldots,\pi^{m}$  the local inverses of $\pi$. Then
$$\overline{J\pi(w)}K^{\phi}_{G}(\pi(z),\pi(w))=\sum_{k=1}^{m}K^{\phi\circ \pi}_{D}(\pi^{k}\circ \pi(z),w) J\pi^{k}(\pi(z))$$
for any $z\notin \pi^{-1}(\pi(N(J\pi)))$, where $ N(J\pi)=\{J\pi=0\}$.
\end{corollary}

The idea of the definition of the operator $\Gamma$ introduced above comes from \cite{Misra} where the Authors discussed the case of the symmetrization mapping between the polydisc and the symmetrized polydisc. In our paper we show that those methods work in general and in particular we show in details how we may apply this method to another domain-tetrablock that is defined in a similar manner to the symmetrized polydisc. 
Recall the definition of that object. Let $\varphi:\mathcal{R}_{II}\rightarrow \mathbb{C}^{3}$, $\varphi(z_{11},z_{22},z):=(z_{11},z_{22},z_{11}z_{22}-z^{2})$ where $\mathcal{R}_{II}$ denotes the Cartan domain of the second type, that is the set of all symmetric matrices of order 2 with the operator norm smaller than $1$ (we identify $(z_{11},z_{22},z)$ with $2\times 2$ symmetric matrix 
$ \left( \begin{array}{ll}
                 z_{11}   &  z \\
                 z   &  z_{22}
\end{array} \right) ). $ 
Then  $\varphi$ is a proper holomorphic map and $\varphi(\mathcal{R}_{II})=\mathbb{E}$ is a domain, called {\it tetrablock}. We find the effective formula for the Bergman kernel of tetrablock and show that tetrablock is not the Lu Qi-Keng domain. Recall that the domain $D$ is {\it a Lu Qi-Keng domain} if its Bergman kernel does not have zeros and is not a Lu Qi-Keng domain if it has.  

\begin{corollary} For any $z,\,w\in \mathbb{R}_{II}$
 $$J\varphi(z)\overline{J\varphi(w)}\,K_{\mathbb{E}}(\varphi(z),\varphi(w))= 
\frac{1}{2}\Big{(}K_{\mathcal{R}_{II}}\Big{(}(z_{11},z_{22},z),w\Big{)}-K_{\mathcal{R}_{II}}\Big{(}z_{11},z_{22},-z),w\Big{)}\Big{)}.$$
\end{corollary}
A consequence of the last formula is the following.
\begin{corollary}
 $\mathbb{E}$ is not a Lu Qi-Keng domain.
\end{corollary}
One may show directly the fact that $\varphi$ described above is a proper holomorphic mapping and $\mathbb{E}$ is a domain. However, it seems reasonable to formulate a result which will be the generalization of that fact and that will help us to avoid the ad hoc proof of properness and opennes of a wide class of mappings. That is the reason why we present below some auxiliary result whose idea of the proof basically comes from the proof of Proposition 2.1 in \cite{Rudin}. We formulate and show it in a more general setting so that it could be applied among others to the above mentioned case of the tetrablock and the symmetrized polydisc.
\begin{propozycja}
Let $\pi:D\rightarrow \mathbb{C}^{n}$ be a holomorphic map where $D\subset\mathbb{C}^{n}$ is a domain and $\mathcal{U}$ is a finite group of homeomorphic transformations of $D$ such that $D$ is precisely $\mathcal{U}$-invariant, that is for $z,\,w\in\mathbb{C}^{n}$ we get that $\pi(z)=\pi(w)$ if and only if $Uz=w$ for some $U\in\mathcal{U}$. Then $\pi(D)$ is a domain and $\pi:D\rightarrow \pi(D)$ is a proper mapping.
\end{propozycja} 
\begin{uwaga}
In Proposition 1 we only assumed that every $U\in\mathcal{U}$ is a homeomorphism but the equality $\pi\circ U=\pi$ easily implies that actually $\mathcal{U}$ is necessarily contained in the group of holomorphic automorphisms of $D$.
\end{uwaga}
\begin{uwaga}
Map $\varphi$ (defined above) is $\mathcal{U}_{\mathbb{E}}:=\{\textup{Id},\,\textup{diag}(1,1,-1)\}$-invariant. We need to verify whether $\mathcal{U_{\mathbb{E}}}$ describes subgroup of group of automorphisms of $\mathcal{R}_{II}$. It can be derived from two facts: 1)every symmetric matrix can be represented in the form $V\Lambda V^{t}$, where $V$ is unitary and $\Lambda=\textup{diag}(\lambda,\mu)$ with $\lambda\geq \mu\geq 0$ (see \cite{Hua} p. 63); and 2) $\lVert A\rVert^{2}=r(AA^{*})$, where $r$ is the spectral radius. Note that the formula for eigenvalues of $AA^{*}$, where $A$ is symmetric, involves only $z_{11},\,z_{22},\,z^{2}$.
\end{uwaga}
\begin{uwaga}
Let us consider $\pi=(\pi_{1},\ldots,\pi_{n})$, where $\pi_{j}$ is the $j$-th elementary symmetric polynomial. In that case the finite group of unitary transformations under which $\pi$ is precisely invariant is the group of permutations $\mathcal{S}_{n}$. Proposition 1 gives the proof of the fact that $\pi|_{\mathbb{D}^{n}}$ is a proper holomorphic mapping onto the image i. e. the symmetrized polydisc and the symmetrized polydisc is open. 
\end{uwaga}
\begin{uwaga}
Fix any $k>2$ and consider: function $\varphi_{k}:\mathcal{R}_{II}\rightarrow \mathbb{C}^{3}$, $\varphi_{k}(z_{11},z_{22},z)=(z_{11},z_{22},z_{11}z_{22}-z^{k})$ and set $\varphi_{k}(\mathcal{R}_{II})$. Notice that $\varphi$ is not proper onto its image. If it were, then the map $\mathcal{R}_{II}\rightarrow \mathcal{R}_{II}$,  $(z_{1},z_{2},z)\rightarrow (z_{1},z_{2},\zeta z)$ should be an isometry (with respect to
the operator norm) for every $\zeta^{k}=1$ but simple examples show that the last one does not hold. 

We might go further and consider Cartan domain of second type $\mathcal{R}_{II}$ in $\mathbb{C}^{\binom{n}{2}}$, here $\mathcal{R}_{II}$ is the set of all symmetric matrices of order n with the operator norm smaller than $1$ (for definitions and properties see \cite{Hua}), and investigate holomorphic map $\varphi:\mathcal{R}_{II}\rightarrow \mathbb{C}^{\binom{n+1}{2}},\ 
\varphi((z_{jk})_{1\leq j\leq k\leq n})=(z_{1,1},\ldots,z_{n,n},z_{1,1}z_{2,2}-z^{2}_{1,2},z_{1,1}z_{3,3}-z^{2}_{1,3},\ldots,z_{1,1}z_{n,n}-z^{2}_{1,n},\ldots,z_{n-1,n-1}z_{n,n}-z^{2}_{n-1,n})$. Unfortunatelly, $\varphi$
fails to be proper onto the image (for the same reason as $\varphi_{k}$ are not), either. So, this indicates that there is no obvious generalization of the tetrablock in higher dimension.
\end{uwaga}

The tetrablock was first studied in \cite{A}. Afterwards it was studied by many authors. In particular, it was shown that the tetrablock is a $\mathbb{C}$-convex domain (see \cite{Zwonek}). The importance of the tetrablock for the geometric function theory follows from the fact that it is the second example (the first one was was the symmetrized bidisc) which is hyperconvex and not biholomorphically equivalent to a convex domain but despite it the Lempert Theorem (see \cite{Lempert1} and \cite{Lempert2}) holds for it (see \cite{Edigarian}). It is also natural to find the Bergman kernel for the tetrablock (using the formula for the Bergman kernel of the Cartan domain and Bell's transformation formula). To our surprise it turned out the tetrablock is not the Lu Qi Keng domain; moreover, it vanishes at very simple points. 

As to the history of the Lu Qi Keng problem we refer the interested Reader to \cite{Boas}. There are many results on  both : domains being Lu Qi-Keng and being not Lu Qi-Keng (see e.g. in \cite{Boas2}, \cite{Zhang}) . 

Recall that (\cite{Hua} p. 84)
$$K_{\mathcal{R}_{II}}(t,s)=\frac{1}{\textup{Vol}(\mathcal{R}_{II})}\Big{(}\textnormal{det}(I-t\overline{s})\Big{)}^{-3}$$
for $t,\,s\in\mathcal{R}_{II}$. Since every point in $\mathcal{R}_{II}$ can be carried by some automorphism of $\mathcal{R}_{II}$ into origin (see \cite{Hua} p. 84), we get $K_{\mathcal{R}_{II}}\not= 0$. Thus $\mathcal{R}_{II}$ is a Lu Qi-Keng domain. Therefore, we have a proper holomorphic mapping $\varphi:\mathcal{R}_{II}\rightarrow\mathbb{E}$ of multiplicity 2 such that $\mathcal{R}_{II}$ is a Lu Qi-Keng domain whereeas $\mathbb{E}$ is not Lu Qi-Keng domain. Recall that another example of that type is $\{|z|+|w|<1\}\owns (z,w)\rightarrow (z^{2},w)\in\{|z|^{\frac{1}{2}}+|w|<1\}$ (see \cite{Boas2}). In our situation there is equality of holomorphically invariant distances in both domains and both domains are $\mathbb{C}$-convex (see \cite{Edigarian}, \cite{Zwonek}) whereas in the example from \cite{Boas2} it is not the case.

\subsection{Proof of main result}

\begin{proof}[Proof of Proposition 1]
The idea of the proof is based on \cite{Misra}. Recall that the operator 
$$\Gamma:\mathbb{A}_{\phi}^{2}(G)\rightarrow \mathbb{A}_{\phi \circ \pi}^{2}(D)$$
is defined as follows
$$\Gamma(f)(z)=\frac{1}{\sqrt{m}}(f\circ \pi)(z)J\pi(z)$$ 
where $f\in\mathbb{A}_{\phi}^{2}(G)$ and $z\in D$. The idea of the formula of $\Gamma$ comes from the formula
\begin{equation}m\int_{G}f\phi\,dV=\int_{D}(f\circ\pi)|J\pi|^{2}(\phi\circ \pi)\,dV\ \ \textup{for any }f\in \textup{L}_{\phi}^{2}(G). \label{eq:zmiana}
\end{equation}
which makes $\Gamma$ an isometry, so the range of $\Gamma$ is a closed Hilbert subspace of $\mathbb{A}_{\phi\circ \pi}^{2}(D)$. Therefore $\Gamma$ is a unitary operator from $\mathbb{A}^{2}_{\phi}(G)$ onto $\Gamma\mathbb{A}_{\phi}^{2}(G)$. 

Thus, there is an orthogonal projection $P$ from $\mathbb{A}_{\phi \circ\pi}^{2}(D)$ onto $\Gamma\mathbb{A}_{\phi}^{2}(G)$. We prove that $P$ is given by the formula 
$$Pg=\frac{1}{m}\sum_{k=1}^{m}(g\circ \pi^{k}\circ\pi)J(\pi^{k}\circ \pi)$$
where $g\in\mathbb{A}_{\phi\circ\pi}^{2}(D)$. Let us denote the right side by $\widetilde{P}$. First, we need to show that $\widetilde{P}$ is well defined. For this purpose, using the properness of $\pi$ one can easily compute  
\begin{multline*}
\lVert \widetilde{P}g\rVert^{2}_{\mathbb{A}_{\phi\circ\pi}^{2}(D)}=\frac{1}{m^{2}}\int_{D}\Big{|}\sum_{k=1}^{m}(g\circ\pi^{k}\circ\pi)J(\pi^{k}\circ\pi)\Big{|}^{2}(\phi\circ\pi)dV\leq  \\
\frac{1}{m}\int_{D}\sum_{k=1}^{m}|(g\circ\pi^{k}\circ\pi)J(\pi^{k}\circ\pi)|^{2}(\phi\circ\pi)dV=\lVert g\rVert^{2}_{\mathbb{A}_{\phi\circ\pi}^{2}(D)}
\end{multline*}
for $g\in\mathbb{A}_{\phi\circ\pi}^{2}(D)$. It remains to verify whether $\widetilde{P}g$ is holomorphic. For this notice the map $\frac{\widetilde{P}g}{J\pi}$ is a well defined holomorphic function on the set $D\setminus \pi^{-1}(\pi(N(J\pi)))$, constant on the fibres of $\pi$, that is on the set $\pi^{-1}(\pi(z))$ for any $z\in D$. So it induces some map $\widetilde{(\frac{\widetilde{P}g}{J\pi})}$ which is holomorphic on $G \setminus \pi(N(J\pi))$. Riemann Singularity Theorem (see. e.g. \cite{Jarnicki} p. 369) finishes the correctness of the definition of $\widetilde{P}$ provided  we know that $\widetilde{(\frac{\widetilde{P}g}{J\pi})}$ is square integrable on $G$ with weight $\phi$. But for that it is enough to show that $\widetilde{P}g\in\mathbb{A}_{\phi\circ\pi}^{2}(D)$ what we have just proved. Note that we have proved that  
for any $g\in\mathbb{A}_{\phi\circ\pi}^{2}(D)$ the equation $\widetilde{P}g=\Gamma f$ has solution $f$ in $\mathbb{A}_{\phi}^{2}(G)$.
Secondly, notice that $\widetilde{P}^{2}=\widetilde{P}$. Indeed, 
\begin{multline*}
\widetilde{P}^{2}g=\frac{1}{m}\sum_{l=1}^{m}(\widetilde{P}g\circ \pi^{l} \circ \pi)J(\pi^{l} \circ \pi)=
\frac{1}{m^{2}}\sum_{k,l=1}^{m}(g\circ \pi^{l}\circ \pi\circ \pi^{k}\circ\pi)[J(\pi^{l} \circ \pi)\circ \pi^{k} \circ \pi]J(\pi^{k} \circ \pi)  \\
=\frac{1}{m^{2}}\sum_{k,l=1}^{m}(g \circ \pi^{l}\circ \pi)J(\pi^{l} \circ \pi \circ \pi^{k} \circ \pi)=\frac{1}{m^{2}}\sum_{k,l=1}^{m}(g \circ \pi^{l}\circ \pi)J(\pi^{l} \circ \pi)=\widetilde{P}g 
\end{multline*}
for $g\in\mathbb{A}_{\phi\circ\pi}^{2}(D)$. Up to this point, we know that $\widetilde{P}$ is the projection.  
Next we proceed to show the equality $\textup{ran}\,\Gamma=\textup{ran}\,\widetilde{P}$. Similarly as above we get $\widetilde{P}\circ\Gamma =\Gamma $, which gives "$\subset$".  It remains to prove the oposite inclusion.  So the question is whether $\widetilde{P}$ takes values in $\Gamma\mathbb{A}_{\phi}^{2}(G)$. Since $\widetilde{P}^{2}=\widetilde{P}$ it is enough to show that for any $g\in\mathbb{A}_{\phi\circ\pi}^{2}(D)$ the equation $\widetilde{P}g=\Gamma f$ has solution $f$ in $\mathbb{A}_{\phi}^{2}(G)$, and this is true as we proved it before.
\end{proof}
\begin{proof}[Proof of Corollary 1.]
We keep the notation from the previous proof. 
Recall that the adjoint operator to $\Gamma$ is given by the formula $\Gamma^{*}(g)=\sqrt{m}\,\widetilde{(\frac{g}{J\pi})}$ for $g\in\Gamma^{*}\mathbb{A}_{\phi}^{2}(G)$. To finish it suffices to note that for any $f\in\mathbb{A}_{\phi}^{2}(G)$ and $w\in D$ the following equalities hold
\begin{multline*}
\langle \Gamma f,PK^{\phi\circ\pi}_{D}(\cdot,w)\rangle_{\mathbb{A}_{\phi\circ\pi}^{2}(D)}=\langle \Gamma f,K^{\phi\circ\pi}_{D}(\cdot,w)\rangle_{\mathbb{A}_{\phi\circ\pi}^{2}(D)}=\Gamma f(w)=\frac{1}{\sqrt{m}} f(\pi(w))J\pi(w)  \\
=\langle f,K^{\phi}_{G}(\cdot,\pi(w))\rangle_{\mathbb{A}_{\phi}^{2}(G)}\frac{J\pi(w)}{\sqrt{m}}    =\langle\Gamma f,\Gamma K^{\phi}_{G}(\cdot,\pi(w))\rangle_{\mathbb{A}_{\phi\circ\pi}^{2}(D)}  \frac{J\pi(w)}{\sqrt{m}}.
\end{multline*}
So from the Riesz Representation Theorem (uniqueness) we get
$$\overline{J\pi (w)}K_{G}^{\phi}(\pi(\cdot),\pi(w))=(\Gamma^{*}\circ P)K^{\phi\circ\pi}_{D}(\cdot,w)(\pi(\cdot))$$
and this equality holds on $D\setminus \pi^{-1}(\pi(N(J\pi)))$ for arbitrary $w\in D$. If we take $w\notin \pi(N(J\pi))$ then on the same set we have  
$$K_{G}^{\phi}(\pi(\cdot),\pi(w))=(\Gamma^{*}\circ P)\frac{1}{\overline{J\pi (w)}}K^{\phi\circ\pi}_{D}(\cdot,w).$$ 
\end{proof}

\subsection{Application to tetrablock}
\begin{proof}[Proof of Proposition 1.]
Let $\{K_{k}\}_{k\in\mathbb{N}}$ be an increasing sequence of compact subdomains of $D$ exhausting it. Consider the new sequence $\{D_{k}:=\bigcup_{U\in\mathcal{U}}U(K_{k})\}_{k}$. Note that there is some $N$ such that for $k>N$ the set $D_{k}$ is a relatively compact domain in $D$ and the sequence $\{D_{k}\}_{k}$ is exhausting $D$. Fix $k>N$. Then $\{U|_{D_{k}}: U\in\mathcal{U}\}$ is a finite group of automorphisms of $D_{k}$ and $\pi |_{D_{k}}:D_{k}\rightarrow \mathbb{C}^{n}$ is precisely $\mathcal{U}$-invariant. These two facts together with the properness of $U|_{D_{k}}$ as a selfmap of $D_{k}$ imply that the set $\pi(D_{k}) \cap \pi(\partial D_{k})$ is empty.
Let $\Omega_{k}$ be the component of $\mathbb{C}^{n}\setminus\pi(\partial D_{k})$  that contains $\pi(D_{k})\subset \Omega_{k}$ and  consequently $\pi(\partial D_{k})\subset\partial{\Omega_{k}}$. This implies that $\pi|_{D_{k}}:D_{k}\rightarrow\Omega_{k}$ is a proper map.

Therefore, $\Omega_{k}=\pi (D_{k})$. Let $\Omega=\cup_{k}\Omega_{k}$ then $\Omega=\pi(D)$ is a domain in $\mathbb{C}^{n}$. The properness of $\pi$ might be checked as follows. If $K\subset \Omega$ is compact then $K\subset \Omega_{k}$ for some $k$. Hence $\pi^{-1}(K)$ is a compact subset of $D_{k}$ and thus a compact subset of $\Omega$.
\end{proof} 

In some special cases it is easy to describe explicity subspace $\Gamma\mathbb{A}^{2}(G)$ looking at the Taylor expansion. For instance when $\pi$ is the symmetrization map from the unit polydisc to the symmetrized polydisc then $\Gamma\mathbb{A}^{2}(\mathbb{G}_{n})$ consists of all antysymmetric functions from $\mathbb{A}^{2}(\mathbb{D}^{n})$ (for details s.e. \cite{Misra}). Another example (the case of the tetrablock) is contained in the next proof, where we obtain formula for $K_{\mathbb{E}}$, which might also be obtained explicitly from Bell's formula. 

\begin{proof}[Proof of Corollary 2.]
We keep the notation from the proof of Theorem 1.
The range of the operator $\Gamma$ contains precisely those maps whose coefficients at $z_{11}^{k}z_{22}^{l}z^{2n}$ in the Taylor expansion 
at the origin vanish for $k,\, l,\, n$ natural numbers (see below). We showed that every function in $\Gamma\mathbb{A}^{2}(\mathbb{E})$ is of the form $J\pi\cdot h$ for some function $h$ depending on $z_{11},\,z_{22},\,z^{2}$. Projection
$$P:\mathbb{A}^{2}(\mathcal{R}_{II})\rightarrow\mathbb{A}^{2}_{z}(\mathcal{R}_{II})$$
acts as follows
$$P(f)(z_{11},z_{22},z)=\frac{1}{2}\Big{(}f(z_{11},z_{22},z)-f(z_{11},z_{22},-z)\Big{)}.$$
Finally,
\begin{multline*}
K_{\mathbb{E}}(\varphi(z_{11},z_{22},z),\varphi(w_{11},w_{22},w))=    \\ \frac{K_{\mathcal{R}_{II}}\Big{(}(z_{11},z_{22},z),(w_{11},w_{22},w)\Big{)}-K_{\mathcal{R}_{II}}\Big{(}(z_{11},z_{22},-z),(w_{11},w_{22},w)\Big{)}}{2J\varphi(z_{11},z_{22},z)\overline{J\varphi(w_{11},w_{22},w)}}
\end{multline*}
for $(z_{11},z_{22},z),\,(w_{11},w_{22},w) \notin N(J\varphi)$.
\end{proof}

\begin{proof}[Proof of Corollary 3.]
We examine a formula for the Bergman kernel for pair $\varphi(0,0,1),\,\varphi(0,0,z)$ (note that the formula for the Bergman kernel extends analytically to $\overline{\mathcal{R}_{II}}\times\mathcal{R}_{II}$). Easy calculation shows that $\frac{-12}{\pi^{3}}K_{\mathbb{E}}(\varphi(0,0,1),\varphi(0,0,z))=6+20z^{2}+6z^{4},\ z\in\mathbb{D}$. The last expresion has zero $z_{0}\in\mathbb{D}$ (s.t. $z_{0}^{2}\in (-1,0)$). 
Now  the equality $K_{\mathbb{E}}(\varphi(0,0,1),\varphi(0,0,z_{0}))=K_{\mathbb{E}}(\varphi(0,0,r),\varphi(0,0,\frac{1}{r} z_{0}))$ for $0<r<1$ s.t. $\frac{z_{0}}{r}\in\mathbb{D}$ finishes the proof. 
\end{proof}

\end{document}